\documentclass[11pt]{article} 

\usepackage{a4,amssymb, amsmath, amsthm}
\usepackage{color}
\usepackage[utf8]{inputenc}
\usepackage[T1]{fontenc}
\usepackage{mathtools}
\usepackage{nicefrac}

\usepackage[a4paper, total={6in, 8in}]{geometry}
 
 \usepackage{graphicx,color}
\usepackage[colorlinks=true]{hyperref}
 \usepackage{csquotes}
 
 \usepackage{amsmath,amssymb,amsthm,mathrsfs,amsfonts,dsfont}
 
\newtheorem{thm}{Theorem}[section]

\newtheorem{lem}[thm]{Lemma}
\newtheorem{pro}[thm]{Proposition}

\newtheorem{defn}[thm]{Definition}

\newcommand {\bv} { {\mathbf{v}} }
\newcommand {\bx} { {\mathbf{x}} }
\newcommand{\diff}{\textnormal{d}}

\theoremstyle{remark}
\newtheorem{rem}[thm]{Remark}

\usepackage{graphicx,epstopdf}
\ifpdf
  \DeclareGraphicsExtensions{.eps,.pdf,.png,.jpg}
\else
  \DeclareGraphicsExtensions{.eps}
\fi

\usepackage{subfig}

\begin{document}

\title{Macroscopic boundary conditions for a fractional diffusion equation in chemotaxis\\\vskip 0.8cm}
\author{Gissell Estrada-Rodriguez\thanks{Universitat Politecnica de Catalunya, Barcelona, 08034, Spain. Centre de Recerca Matematica, 08193 Bellaterra, Barcelona, Spain} \and  Heiko Gimperlein\thanks{Engineering Mathematics, University of Innsbruck, 6020 Innsbruck,  Austria. \newline \noindent email: gissell.estrada@upc.edu, heiko.gimperlein@uibk.ac.at.\newline \noindent G.~E.~R.~was supported by the Spanish grant PID2022-143012NA-I00 and she is member of the Catalan research group 2021-SGR-00087. This work is supported by the Spanish State Research Agency, through the Severo Ochoa and María de Maeztu Program for Centers and Units of Excellence in R\&D (CEX2020-001084-M)} }
\date{}

\maketitle \vskip 0.5cm
\begin{abstract}
\noindent In this paper we examine boundary effects in a fractional chemotactic equation derived from a kinetic transport model describing cell movement in response to chemical gradients (chemotaxis). Specifically, we analyze reflecting boundary conditions within a nonlocal fractional framework. Using boundary layer methods and perturbation theory, we derive first-order approximations for interior and boundary layer solutions under symmetric reflection conditions. This work provides fundamental insights into the complex interplay between fractional dynamics, chemotactic transport phenomena, and boundary interactions, opening  future research in biological and physical applications involving nonlocal processes.  
\end{abstract}

\vskip 1.0cm

\section{Introduction}

Fractional partial differential equations extend classical models by incorporating nonlocal movement through fractional derivatives, which modify the formulation and impact of boundary conditions. In contrast to traditional local operators where boundary effects are confined to regions immediately adjacent to the boundary, fractional operators require the solution to be defined both within and outside the domain. This nonlocal nature means that mass situated deep within the domain can still be influenced by the boundary, leading to rich and diverse phenomena not observed in classical settings.

Recent numerical investigations have highlighted these differences in the context of space-fractional diffusion equations on the unit interval. For example, Baeumer et al.~\cite{baeumer2018reprint} examined both absorbing (Dirichlet) and reflecting (Neumann) boundary conditions. Their study showed that when the exponent of the fractional derivative 
$\alpha=2$, the reflecting boundary condition converges to the classical Neumann condition. In another work, Baeumer et al. \cite{baeumer2018fractional} introduced a “fast-forward” reflective boundary condition, where the time that mass spends outside the domain is effectively neglected, allowing for an accelerated re-entry into the domain. These approaches highlight the necessity of revisiting traditional boundary condition concepts when dealing with fractional operators. 
In a more recent work, You et al. \cite{you2020asymptotically} introduced an asymptotically compatible method for imposing Neumann-type boundary conditions on two-dimensional nonlocal diffusion problems, offering a new formulation that achieved optimal second-order convergence to the local limit as the nonlocal horizon parameter approaches zero.
Their method carefully constructed a nonlocal flux boundary condition that is rigorously shown to converge to the classical Neumann problem, ensuring well-posedness and avoiding artificial surface effects, even in non-convex domains.

In this paper, we focus on the fractional chemotactic equation \cite{pks}
\begin{equation}
\partial_tu=\nabla\cdot(C_\alpha\nabla^{\alpha-1}u-\chi u\nabla\rho)\ ,\label{eq: main at intro}
\end{equation}
where $\nabla^{\alpha-1}$ denotes a fractional gradient which interpolates between ballistic motion ($\alpha=1$) and ordinary diffusion ($\alpha=2$); note that the case $\alpha=2$ corresponds to the classic formulation of Patlak--Keller--Segel equations proposed phenomenologically in \cite{keller1970initiation}. In \eqref{eq: main at intro}, the chemotactic population is governed by a diffusion term with coefficient $C_\alpha$ (defined in \eqref{eq: constants})
that represents a random component to motility, and a chemotactic flux of advective type, 
where the advection is proportional to the chemical gradient. The function 
$\chi$ is commonly referred to as the chemotactic sensitivity. In the case of constant $C_\alpha$ we obtain an honest fractional Laplacian, namely, $\nabla\cdot\nabla^{\alpha-1}=c(-\Delta)^{\nicefrac{\alpha}{2}}$ for $1<\alpha<2$.

 Chemotaxis, the directed movement of cells or organisms in response to chemical gradients, is a phenomenon that becomes particularly complex when nonlocal transport mechanisms are considered. By extending the classical chemotactic model to include fractional dynamics, we capture the subtleties of mass transport that are influenced by long-range movement and boundary effects. We explore the implications of reflecting boundary conditions in this nonlocal framework.

To address these challenges, based on classical works by Alt \cite{alt1980singular}, respectively Larsen \cite{larsen}, our analysis employs boundary layer techniques and perturbation theory to derive first-order approximations for both the interior and boundary layers of the solution starting from a kinetic description of the movement. Let $f(\bx,t,\bv)$ be a probability density function which describes the distribution in space $\bx\in\Omega\subset\mathbb{R}^n$ and velocity $\bv\in S=\{\bv\in\mathbb{R}^n:\ |\bv|=1 \}$ at time $t$ of individuals, then the evolution of $f$ is given by \cite{alt1980singular,robots,pks} 
\begin{align}
\partial_t{f}(\cdot,\mathbf{v})+c\mathbf{\mathbf{v}}\cdot \nabla {f}(\cdot,\mathbf{v}) =-(\mathds{1}-T)\int_0^{t}(\beta \bar{f})(\cdot,\mathbf{v},\tau) \diff\tau\ . \label{eq: 3.15}
\end{align}
Here $\tau\in\mathbb{R}$ describes the run time, i.e. the time of approximately straight line motion before the individual stops and chooses a new direction, and  the turn angle operator $T$ is given by
\begin{equation}
T\phi(\eta) =\int_S k(\mathbf{x},t,\mathbf{\mathbf{v}};\mathbf{\eta})\phi(\mathbf{v})\diff\mathbf{v} = \int_S \ell(\mathbf{x},t,|\eta-\mathbf{v}|)\phi(\mathbf{v})\diff\mathbf{v} \label{eq: turn angle operator}\ ,
\end{equation}
which describes the effect of changing from direction $\bv$ to a new direction $\mathbf{\eta}$.

The investigation is carried out under the assumption of symmetric reflection at the boundary, which, while a simplification, offers valuable insights into the underlying mechanics of the system. This approach not only provides a theoretical framework for understanding the boundary effects in fractional chemotaxis but also bridges the gap between earlier numerical results on space-fractional diffusion and the emerging study of fractional chemotactic phenomena.

 By integrating established numerical findings with our analytical results, we set the stage for further exploration into the diverse applications and implications of fractional models in both biological and physical systems.

\section{Modelling equations and scaling}\label{sec: fractional PKS}

Motivated by the experimental results in \cite{korobkova2004molecular} and \cite{li2008persistent}, we model a population of organisms moving in a medium in $\Omega\subset\mathds{R}^n$, containing 
some chemical (with concentration $\rho=\rho(\mathbf{x},t)$) 
that acts as an attractant. We assume that each individual performs a biased random walk with the following properties:
\begin{enumerate}
\item Starting at position $\mathbf{x}$ and time $t$, we assume an individual runs in 
direction $\mathbf{v}$ for some time $\tau$, called the \enquote{run time}. 
\item The individuals are assumed to move with constant forward speed $c$, following a straight line motion between reorientations. 
\item Each time the individual stops it selects a new direction $\eta$ according to 
a distribution $k(\mathbf{x},t,\mathbf{\mathbf{v}};\eta)$, which only depends on $|\mathbf{v} - \eta|$. The choice 
of new direction is taken here to be independent of the chemical concentration or gradient. 
\item The reorientation is assumed to be (effectively) instantaneous.  
\end{enumerate}

Now assume that $\mathcal{X}$ and $\mathcal{T}$ are the macroscopic space and time scales respectively. Let us also consider that the mean run time $\bar{\tau}$ is small compared with the macroscopic time $\mathcal{T}$, i.e., $\epsilon=\nicefrac{\bar{\tau}}{\mathcal{T}}\ll 1$ where $\epsilon$ is a small parameter. Suppose further that the concentration $\rho$ is already dimensionless in the sense that it stands for $\rho/\rho_0$ where $\rho_0$ is an averaged value of $\rho$ over $\Omega$. 

We consider the scaling $$t_n=\epsilon t\ ,\quad  \mathbf{x}_n=\frac{\epsilon \mathbf{x}}{s}\ ,\quad  c_n=\epsilon^{-\gamma} c_0\quad \textnormal{and}\quad  \tau_n=\tau\epsilon^{\mu}\ ,$$ 
for $\mu>0$ and $0<\gamma<1$.
Equation \eqref{eq: 3.15} becomes,
\begin{equation}
\epsilon\partial_t{f}+\epsilon^{1-\gamma}c_0\mathbf{\mathbf{v}}\cdot \nabla {f}=-(\mathds{1}-T)\int_0^{t}\beta_\epsilon \bar{f} \diff\tau \label{eq: 3.13}\ ,
\end{equation}
where $\bar{f}=\bar{f}(\bx,t,\bv,\tau)$ and 
\begin{equation}
\beta_\epsilon(\cdot,\mathbf{v},\tau)=\frac{\alpha\epsilon^{\mu}}{\tau_0\epsilon^{\mu}+\tau_1\epsilon^{\mu}D^\mathbf{v}_\epsilon\rho+\tau}\ .\label{eq: master with scaling}
\end{equation}
Here the total gradient of $\rho$ is given by $
D^\mathbf{v}_\epsilon\rho=\epsilon\partial_t \rho+\epsilon^{1-\gamma} c_0\mathbf{v}\cdot \nabla\rho.
$

Since $\ell$ in (\ref{eq: turn angle operator}) is a probability distribution, it is normalized to $\int_S \ell(\mathbf{x},t,|\mathbf{v}-e_1|)\diff\mathbf{v}=1$,  where $e_1 = (1,0,\dots, 0)$. We immediately observe
\begin{align}
    \int_S (\mathds{1}-T)\phi \diff\mathbf{v}=0\label{eq: conservation T 2}\ ,
\end{align}
for all $\phi\in L^2(S)$. Biologically, (\ref{eq: conservation T 2}) corresponds to the conservation of the number of individuals in the tumbling phase.

We also require some more detailed information about the spectrum of $T$. Recall that in $n$-dimensions, the surface area of the unit sphere $S$ is given by \[
|S|=\begin{cases}
\frac{2\pi^{\nicefrac{n}{2}}}{\Gamma\left(\frac{n}{2}\right)}\ , & \textnormal{for $n$ even}\ ,\\
\frac{\pi^{\nicefrac{n}{2}}}{\Gamma\left(\frac{n}{2}+1\right)}\ , & \textnormal{for $n$ odd}\ .
\end{cases}
\]
\begin{lem}\label{lem: eigenfunctions}
Assume that $\ell$ is continuous. Then $T$ is a symmetric compact operator. In particular, there exists an orthonormal basis of $L^2(S)$ consisting of eigenfunctions of $T$.\\ With $\mathbf{\mathbf{v}}=(\mathbf{v}_0,\mathbf{v}_1,...,\mathbf{v}_{n-1}) \in S$, we have
\begin{equation}
\begin{aligned}\phi_{0}(\mathbf{v}) & =\frac{1}{|S|} &  & \text{is an eigenfunction to the eigenvalue} &  & \nu_{0}=1\ ,\\
\phi_{1}^j(\mathbf{v}) & =\frac{n\mathbf{v}_j}{|S|} &  & \text{are eigenfunctions to the eigenvalue} &  & \mathbf{u}_{1}=\int_{S}\ell(\cdot,|\mathbf{u}-1|)\mathbf{u}_{1}\diff\mathbf{u}<1\ . \label{eq: eigen}
\end{aligned}
\end{equation}
Any function $\bar{f}\in L^2(\mathds{R}^n\times \mathds{R}^+\times S)$ admits a unique decomposition 
\begin{equation}
{f}=\frac{1}{|S|}\left({u}+n\mathbf{\mathbf{v}}\cdot {w} \right)+\hat{z},\label{eq: real_eigen}
\end{equation}
where $\hat{z}$ is orthogonal to all linear polynomials in $\mathbf{v}$. Explicitly,
\begin{equation*}
{u}(\mathbf{x},t)=\int_S{f}(\mathbf{x},t,\mathbf{\mathbf{v}})\phi_0(\mathbf{v}) \diff\mathbf{v},\ {w}^j(\mathbf{x},t)=\int_S {f}(\mathbf{x},t,\mathbf{\mathbf{v}})\phi_1^j(\mathbf{v}) \diff\mathbf{v},
\end{equation*}
and  ${w} =  ({w}^1, \dots, {w}^n)$.
\end{lem}

Using the Laplace transform and convolution properties, we can re-write the right hand side of \eqref{eq: 3.13} in terms of the density ${f}$ and a turn angle kernel $\mathcal{T_\epsilon}$, as was shown in \cite{pks}.
Therefore we have the following kinetic equation independent of the run time $\tau$,
\begin{align}
    \epsilon\partial_t{f}+\epsilon^{1-\gamma}c_0\mathbf{v} \cdot\nabla{f} =-(\mathds{1}-T)\mathcal{T}_\epsilon{f}\ ,\label{eq: last}
\end{align}
where, to leading order for $\mu>0$,
\begin{equation}
\begin{aligned}
    \mathcal{T}_\epsilon & =\frac{\epsilon^{-\mu}(\alpha-1)}{\tau_0}  -\frac{\tau_1(\alpha-1)}{\tau_0^2}\epsilon^{-\mu-\gamma+1}c_0\mathbf{v}\cdot\nabla\rho-\frac{\epsilon^{1-\gamma}c_0}{2-\alpha}(\mathbf{v}\cdot\nabla)\\  & -\tau_0^{\alpha-2}\epsilon^{\mu(\alpha-2)+(1-\gamma)(\alpha-1)}(\alpha-1)^2\Gamma(-\alpha+1)(c_0\mathbf{v}\cdot\nabla)^{\alpha-1}+\mathcal{O}(\epsilon^{\mu-\gamma+1})\ .\label{eq: almost final}
\end{aligned}
\end{equation}
Here we have used the fact that $D^\mathbf{v}_\epsilon\rho=\epsilon\partial_t\rho+\epsilon^{1-\gamma}c_0\mathbf{v}\cdot\nabla\rho \simeq \epsilon^{1-\gamma}c_0\mathbf{v}\cdot\nabla\rho$ since $1-\gamma<1$.

\section{Inner solution}\label{sec: inner solution}
We are going to study the following initial value problem:
\begin{equation}
\begin{aligned}\epsilon\partial_t{f}+\epsilon^{1-\gamma}c_0\mathbf{v}\cdot\nabla{f}& =-(\mathds{1}-T)\mathcal{T}_\epsilon{f}\ , &  & \mathrm{in}\ \mathds{R}_+\times{\Omega}\times S\ ,\\
{f}(\mathbf{x},0,\mathbf{v}) & =\delta_0(\mathbf{x})\ , &  & \mathrm{in}\ {\Omega}\times S\thinspace,\\
{f}(\xi,t,\mathbf{v}) & ={f}(\xi,t,\textnormal{R}_\xi(\mathbf{v})) \ ,&  & \ \forall\  \xi\in\partial\Omega\ ,\ \mathbf{v}\in S\ , \ {\mathbf{v}\cdot \nu>0}\ ,
\end{aligned}
\label{eq: original system}
\end{equation}
where $\textnormal{R}_\xi(\mathbf{v})=\mathbf{v}-2(\nu_\xi\cdot\mathbf{v})\nu_\xi$ for $\nu_\xi$ the unit inner normal vector at $\xi\in\partial \Omega$, where $\partial\Omega$ is a smooth boundary. We are going to study the behaviour of a smooth solution $f$ to the previous problem.

Let
\begin{equation}
    {f}^i=\sum_{m=0}^N\varepsilon^{m}{f}_{m}^i\ ,
\end{equation}
be the interior approximation of the solution of the system (\ref{eq: original system}) where $\varepsilon=\epsilon^{1/2}$. We expand  ${f}_{m}^i$ using the eigenfunction representation,
\begin{equation}
{f}^i_{m}=\frac{1}{|S|}({u}_{m}^i+n\mathbf{\mathbf{v}}\cdot{w}^i_{m})+\textnormal{l.o.t.}\ ,\label{eq: expansion}
\end{equation}
where 
\begin{equation*}
{u}_m^i(\mathbf{x},t)=\int_S{f}^i_m(\mathbf{x},t,\mathbf{\mathbf{v}}) \diff\mathbf{v}\ ,\qquad {w}_m^i(\mathbf{x},t)=\int_S \mathbf{v}{f}_m^i(\mathbf{x},t,\mathbf{\mathbf{v}}) \diff\mathbf{v}\ .
\end{equation*}
More precisely, we consider 
\[
f^i=\sum_0^m\varepsilon^m(u_m^i+n\mathbf{v}w_m^i)=(u_0^i+n\mathbf{v}\cdot w_0^i)+\varepsilon(u_1^i+n\mathbf{v}\cdot w_1^i)+\textnormal{l.o.t.}\ .
\]
Hence we can write the system (\ref{eq: original system}) in terms of the new quantities ${u}_m^i$ and ${w}_m^i$ taking the $\varepsilon^m$-component as follows,
\begin{align}
    \varepsilon^2\partial_t& ({u}^i_{m}  + n\mathbf{v}\cdot{w}_{m}^i)  +c_0\varepsilon^{2-2\gamma}\mathbf{v}\cdot\nabla({u}^i_{m}+ n\mathbf{v}\cdot{w}_{m}^i)\nonumber\\ & =-(\mathds{1}-T)\left[\varepsilon^{-2\mu}B_{2\mu}({u}^i_{m}+n\mathbf{v}\cdot{w}^i_{m})+\varepsilon^{1-2\mu}B_{2\mu-1}({u}^i_{m}+n\mathbf{v}\cdot{w}^i_{m})\right].\label{eq: expansion initial}
\end{align}
Here we have rewritten (\ref{eq: almost final}) as 
\begin{equation}
\mathcal{T}_\varepsilon=\varepsilon^{-2\mu}B_{2\mu}+\varepsilon^{-2\mu+1}B_{2\mu-1}, \label{eq: new T}
\end{equation}
where 
\begin{align}
B_{2\mu}&=\frac{(\alpha-1)}{\tau_0}\ ,\label{eq: B0}\\
B_{2\mu-1}&=-\frac{\tau_1(\alpha-1)}{\tau_0^2}c_0\mathbf{v}\cdot\nabla\rho-\tau_0^{\alpha-2}(\alpha-1)^2\Gamma(-\alpha+1)(c_0\mathbf{v}\cdot\nabla)^{\alpha-1}.
\end{align}

For simplicity and following the results in \cite{pks}, we are taking $\gamma=1/2$ and $\mu=\nicefrac{(2-\alpha)}{2(\alpha-1)}$. With this parameter choice, one can check that in (\ref{eq: new T})
\[
\epsilon^{\mu(\alpha-2)+(1-\gamma)(\alpha-1)}=\varepsilon^{-2\mu+1}\ .
\]
By integrating (\ref{eq: expansion initial}) with respect to all directions $\mathbf{\mathbf{v}}$ and using (\ref{eq: conservation T 2}), we find the macroscopic conservation equation
\begin{equation}
\varepsilon^2\partial_t{u}_{m}^i+\varepsilon c_0 n\nabla\cdot{w}_{m}^i=0\ .\label{eq: conservation}
\end{equation}
Now we need to find an expression for the mean direction ${w}_{m}^i$. 
Then we multiply \eqref{eq: expansion initial} by $\mathbf{v}$ and integrate over the unit sphere $S$ to obtain
\begin{align}
\varepsilon^2n\partial_t{w}^i_{m}+\varepsilon c_0\nabla {u}^i_{m} & =-\int_S\mathbf{v}(\mathds{1}-T)\Bigl[\varepsilon^{-2\mu}B_{2\mu}({u}^i_{m}+n\mathbf{v}\cdot{w}^i_{m})\nonumber\\ &\quad +\varepsilon^{1-2\mu}B_{2\mu-1}({u}^i_{m}+n\mathbf{v}\cdot{w}^i_{m})\Bigr]\diff\mathbf{v}.
\end{align}
For $m=0$ we obtain the first non-zero elements of the expansion. Hence, it follows from
\begin{align}
   \varepsilon^{-2\mu}:\hspace{0.5cm}  0=\frac{1}{|S|}\int_S\mathbf{v}(\mathds{1}-T)\Bigl(\frac{\alpha-1}{\tau_0}({u}^i_{0}+n\mathbf{v}\cdot{w}_{0}^i) \Bigr)\diff\mathbf{v},
\end{align}
that ${u}_{0}^i={u}_{0}^i(x)$ and ${w}_{0}^i=0$. From the next power of $\varepsilon$, corresponding to $m=1$ it follows that
\begin{equation}
\begin{aligned}
    \varepsilon^{1-2\mu}:\hspace{0.5cm} &0= -\int_S\mathbf{v}(\mathds{1}-T)\Bigl[\frac{\alpha-1}{\tau_0}({u}_{1}^i+n\mathbf{v}\cdot{w}_{1}^i)+\Bigl(-\frac{\tau_1(\alpha-1)}{\tau_0^2}c_0\mathbf{v}\cdot\nabla\rho\nonumber\\ &\quad -\tau_0^{\alpha-2}(\alpha-1)^2\Gamma(-\alpha+1)(c_0\mathbf{v}\cdot\nabla)^{\alpha-1} \Bigr)({u}_{0}^i+n\mathbf{v}\cdot{w}_{0}^i)\Bigr]\diff\mathbf{v}\ .
\end{aligned}
\end{equation}
Hence we obtain that the flux is given by
\begin{equation}
{w}_{1}^i=\frac{\tau_1}{n\tau_0}c_0{u}_{0}^i\nabla\rho+\frac{\pi\tau_0^{\alpha-1}(\alpha-1)}{\sin(\pi\alpha)\Gamma(\alpha)}\frac{(n^2\nu_1-|S|)}{n|S|(\nu_1-1)}c_0^{\alpha-1}\nabla^{\alpha-1}{u}_{0}^i\ ,\label{eq: mean direction}
\end{equation}
where we have used $\Gamma(-\alpha+1)=\frac{\pi}{\sin(\pi\alpha)\Gamma(\alpha)}$. From the conservation equation (\ref{eq: conservation}) we obtain
\begin{equation}
\partial_t{u}_{0}^i=nc_0\nabla\cdot\left(C_\alpha\nabla^{\alpha-1}{u}_{0}^i-\chi{u}_{0}^i\nabla\rho\right)\label{eq:final equation}\ ,
\end{equation} 
 where
\begin{equation}
    C_\alpha=-\frac{\pi(\tau_0c_0)^{\alpha-1}(\alpha-1)}{\sin(\pi\alpha)\Gamma(\alpha)}\frac{(n^2\nu_1-|S|)}{n|S|(\nu_1-1)}\quad \textnormal{and}\quad \chi=\frac{\tau_1c_0}{n\tau_0}\ .\label{eq: constants}
\end{equation}
We summarize the previous result in the following proposition.
\begin{pro}\label{prop: 1}
Assume $\Omega\subset\mathds{R}^n$ is a bounded domain with smooth boundary $\partial\Omega$. Let ${u}_m^i:\mathds{R}_+\times{\Omega}\rightarrow\mathds{R}$, $0\leq m\leq N$, be the solution of (\ref{eq:final equation}) with initial condition
\[
{u}_{0}^i(\mathbf{x},0)=\delta_0(\mathbf{x})\quad on\quad {\Omega}\ ,
\]
and Neumann (i.e.~reflective) boundary conditions 
\begin{equation}
nc_0\chi\nu\cdot({u}_{0}^i\nabla\rho)+nc_0 \partial_\nu^{\alpha-1} (C_\alpha{u}_{0}^i)=0\quad \textit{on}\quad \mathds{R}_+\times\partial\Omega\ .
\end{equation}
Then, the interior approximation ${f}^i_{(N)}=\frac{1}{|S|}\sum_{m=0}^N\varepsilon^m{f}_m^i$ satisfies the equation (\ref{eq: original system}).

\end {pro}

\section{Half-space boundary approximation}\label{sec: boundary}

In this section we study the system (\ref{eq: original system}) at the boundary. In order to do so, we introduce a correcting function ${f}_m^b$ near the boundary.

Define an asymptotic approximation of the solution at the boundary as ${f}_m^b:\ \mathbb{R}^2_+\times\partial\Omega\times S\rightarrow\mathbb{R}$ for $0\leq m\leq N$ where, for $\xi\in\partial\Omega$,
\begin{equation}
{f}^b=\sum_{m=0}^N\varepsilon^m{f}_m^b(r,\xi,\mathbf{v},t)\ .\label{eq: approximation}
\end{equation}

\begin{defn}[Half-space boundary layer]Since $\partial\Omega$ is compact and smooth, there exists a boundary strip $\Omega_{2\delta}=\left\{x<2\delta \right\}$ such that $r=\textnormal{dist}(x,\partial\Omega)= x/\varepsilon^{\varrho}$ is well defined.
\end{defn}
Note that, in the new coordinates, the gradient is expressed as 
{\[
\nabla =\frac{\partial r}{\partial x}\frac{\partial}{\partial r}+\frac{\partial \xi}{\partial x}\frac{\partial}{\partial \xi} =\frac{1}{\varepsilon^\varrho}\nu\partial_r+\nabla_\xi\ .
\]}
Hence, substituting (\ref{eq: approximation}) into (\ref{eq: last}) we obtain, for $\gamma=1/2$,
\begin{align}
    \varepsilon^2\partial_t{f}^b &+\varepsilon c_0\mathbf{v} \cdot\nabla_\xi{f}^b + \varepsilon^{1-\varrho}c_0\mathbf{v}\cdot\nu\partial_r{f}^b\nonumber\\ & =-\left(\mathds{1}-T \right)\Bigl(\varepsilon^{-2\mu}B_{2\mu}+\varepsilon^{-2\mu+1}A_{2\mu-1} +{\varepsilon^{-2\mu+1-\varrho(\alpha-1)}}C_{2\mu-1+\varrho(\alpha-1)}\Bigr){f}^b\label{eq: big expansion}\ .
\end{align}
In this case $B_{2\mu}$ is given by \eqref{eq: B0} and 
\begin{align}
A_{2\mu-1} & =-\frac{\tau_1(\alpha-1)}{\tau_0^2}c_0\mathbf{v}\cdot\nabla\rho\ ,\ \label{eq: B1}\\
 C_{2\mu-1+\varrho(\alpha-1)}& \simeq-\tau_0^{\alpha-2}(\alpha-1)^2\Gamma(-\alpha+1)(c_0\mathbf{v}\cdot\nu\partial_r)^{\alpha-1}\ ,\label{eq: B2}
\end{align}
where we have used the approximation
{\[
(c_0\mathbf{v}\cdot\nabla)^{\alpha-1}=\varepsilon^{-\varrho(\alpha-1)}(c_0\mathbf{v}\cdot\nu\partial_r)^{\alpha-1}+\mathcal{O}\left(\varepsilon^{-\varrho(\alpha-2)} \right).
\]}
{Now we would like to find $\varrho$ in order to obtain the leading order terms. We notice that in the right hand side of \eqref{eq: big expansion} the leading terms are either the first one (involving $B_{(\cdot)}$) or the third one (involving $C_{(\cdot)}$), depending on $\varrho$. In fact, if we let $\varepsilon^{-2\mu}=\varepsilon^{-2\mu+1-\varrho(\alpha-1)}$ we obtain that $\varrho=2\mu+1=\frac{1}{\alpha-1}$. For $\alpha=2$ (classical diffusion scaling) we get $\mu=0$ and $\varrho=1$, recovering the  scaling used in \cite{alt1980singular} for chemotactic diffusion.
Choosing $\varrho=1/(\alpha-1)$ we write
\begin{align}
    \varepsilon^{2}\partial_t  {f}^b &+\varepsilon c_0\mathbf{v} \cdot\nabla_\xi{f}^b + \varepsilon^{-\frac{2-\alpha}{\alpha-1}}c_0\mathbf{v}\cdot\nu\partial_r{f}^b\nonumber\\ & =-\left(\mathds{1}-T \right)\Bigl(\varepsilon^{-\frac{2-\alpha}{\alpha-1}}B_{\frac{2-\alpha}{\alpha-1}}+\varepsilon^{-\frac{2-\alpha}{\alpha-1}+1}A_{\frac{2-\alpha}{\alpha-1}-1} +\varepsilon^{-\frac{2-\alpha}{\alpha-1}}C_{\frac{2-\alpha}{\alpha-1}}\Bigr){f}^b\ . \label{eq: big expansion final}
\end{align}
Grouping the $\varepsilon$ terms in (\ref{eq: big expansion final}) appropriately we have
 \begin{align*}
     c_0\mathbf{v}\cdot\nu\partial_r  {f}^b+(\mathds{1}-T)\Bigl(B_{\frac{2-\alpha}{\alpha-1}}+C_{\frac{2-\alpha}{\alpha-1}}\Bigr)&{f}^b =\varepsilon^{2\mu}\Bigl(-\varepsilon^2\partial_t{f}^b-\varepsilon c_0\mathbf{v}\cdot\nabla_y{f}^b\\ &-\varepsilon^{1-2\mu}\left(\mathds{1-}T \right)A_{\frac{2-\alpha}{\alpha-1}-1}{f}^b\Bigr)\ ,
 \end{align*}
where $B_{\frac{2-\alpha}{\alpha-1}}$ is a constant given in \eqref{eq: B0}, $C_{\frac{2-\alpha}{\alpha-1}}$ is given in \eqref{eq: B2}  and $A_{\frac{2-\alpha}{\alpha-1}-1}$ is given in \eqref{eq: B1}.
 For convenience of notation, and taking the $\varepsilon^m$-component  we introduce the following result.

 \begin{pro}
     Under the assumptions of this section, and for $\varepsilon$ small, there exist smooth functions $f_m^b:\mathbb{R}^2_+\times\partial\Omega\times S\to\mathbb{R}$ for $1\leq m\leq N$ which satisfy the following system of boundary value problems,
     \begin{equation}
     c_0\mathbf{v}\cdot\nu\partial_r  {f}_{m}^b+(\mathds{1}-T)\Bigl(B_{\frac{2-\alpha}{\alpha-1}}+C_{\frac{2-\alpha}{\alpha-1}}\Bigr){f}^b_{m}=g_{m}^b\ ,\label{eq: m-expression}
 \end{equation}
 where 
 \begin{align}
 g_{m}^b =&-\varepsilon^{2\mu+2}\partial_t{f}_{m}^b-\varepsilon^{2\mu+1}c_0\mathbf{v}\cdot\nabla_\xi{f}^b_{m}-\varepsilon\left(\mathds{1-}T \right)A_{\frac{2-\alpha}{\alpha-1}-1}{f}^b_{m}\ .\label{eq: g}
 \end{align}
 Moreover, we have
 \begin{align}
     {f}_{-2}^b&={f}^b_{-1}=0\ ,\\
     {f}_m^b(0,\xi,\mathbf{v},t) &={f}_m^b(0,\xi,\textnormal{R}_\xi(\mathbf{v}),t) \ .
 \end{align}
 \end{pro}
\iffalse
 Now we aim to find a solution to \eqref{eq: m-expression}. Using Lemma \ref{lem: eigenfunctions} and applying the Fourier transform we write
 \begin{align}
c_0\mathbf{v}\cdot\nu\xi\hat{{f}}_{m+2\mu}^b+(\mathds{1}-\nu_1)B_{2\mu}\hat{{f}}_{m+2\mu}^b=\hat{g}_m^b 
 \end{align}
 \begin{align*}
&\hat{{f}}_{m+2\mu}^b=\frac{\hat{g}_m^b}{c_0\mathbf{v}\cdot\nu\xi+(\mathds{1}-\nu_1)B_{2\mu}}\\
&{{f}}_{m+2\mu}^b={g}_m^b*\frac{1}{c_0\mathbf{v}\cdot\nu}\mathcal{F}^{-1}\Bigl(\frac{1}{\xi+c}\Bigr)
 \end{align*}
 where $c=\frac{(\mathds{1}-\nu_1)B_{2\mu}}{c_0\mathbf{v}\cdot\nu}>0$.
Then we finally obtain
\begin{equation}
{f}^b_{m+2\mu}=\frac{1}{c_0\mathbf{v}\cdot\nu}\int_{-\infty}^\infty g_m^b(r-s)H(s)e^{-cs}d s\ ,
\end{equation}
where $H(r)$ is the Heaviside function.
From \eqref{eq: g} we see that  $\alpha>1$ implies $\alpha(2\mu+1)-2>2\mu-1$ and therefore the term involving $\tilde{B}_{\alpha(2\mu+1)-2}$ dominates over the term involving $B_{2\mu-1}$.
\fi 
While the solution in the interior of the domain satisfies equation (\ref{eq:final equation}), the expansion of the solution at the boundary, ${f}_m^b$, satisfies expression (\ref{eq: m-expression}). \\
 
\noindent To obtain a conservation equation for the new density ${u}_m^b$ at the boundary, we let ${f}_m^b=\frac{1}{|S|}\left({u}_m^b+n\mathbf{v}\cdot{w}_m^b\right)$ in (\ref{eq: m-expression}) and  integrate with respect to all direction in $S$,
\begin{align}\label{eq: int w}
     \frac{c_0}{|S|}\int_S\mathbf{v}\cdot\nu & \partial_r({u}_{m}^b+n\mathbf{v}\cdot{w}_{m}^b)\diff\mathbf{v}\nonumber\\ & =\frac{1}{|S|}\int_S\left(-\varepsilon^{2\mu+2}\partial_t({u}_{m}^b+n\mathbf{v}\cdot{w}_{m}^b)-\varepsilon^{2\mu+1}c_0\mathbf{v}\cdot\nabla_\xi({u}_{m}^b+n\mathbf{v}\cdot{w}_{m}^b) \right)\diff\mathbf{v}.
\end{align}
Thus we obtain
\begin{equation}
-nc_0\partial_r\left(\nu \cdot{w}_{m}^b\right)=\varepsilon^{2\mu+1}nc_0\nabla_\xi\cdot{w}_{m}^b+\varepsilon^{2\mu+2}\partial_t{u}_{m}^b\ .\label{eq: conservation at the boundary}
\end{equation}
Letting $\varepsilon\to 0$ and $m=0$ in \eqref{eq: conservation at the boundary} we observe that $nc_0\partial_r(\nu\cdot w_{0}^b)=0$ which implies that $w_{0}^b$ is constant along the radial direction.
Integrating \eqref{eq: conservation at the boundary} in the boundary strip $\Omega_{2\delta}$, and using 
\[
\frac{1}{\varepsilon^\varrho}\int_{\Omega_{2\delta}}\partial_r\left(\nu\cdot{w}_{m}^b \right)=-\int_{\partial\Omega}\nu\cdot{w}_{m}^b\ ,
\]
from the Fundamental Theorem of Calculus we can write the conservation equation as
\begin{equation}
    \varepsilon^\varrho c_0n\int_{\partial\Omega}\nu\cdot{w}_{m}^b  =\varepsilon^{2\mu+2}\int_{\Omega_{2\delta}}\partial_t{u}_{m}^b\label{eq: mass conservation boundary} \ .
\end{equation}
Expression (\ref{eq: mass conservation boundary}) gives the conservation of mass for the asymptotic expansion at the boundary. \\

Next, we are going to match the interior and boundary solutions. Since we consider that particles are reflected at the boundary, then it holds the necessary condition
\begin{equation}
    \nu\cdot{w}_m^i=-\nu\cdot{w}_m^b\ .\label{eq: continuity at the boundary}
\end{equation}
Integrating the conservation equation (\ref{eq: conservation}) in $\Omega$ and using (\ref{eq: continuity at the boundary}) we get
\begin{equation}
\varepsilon^2\partial_t\int_\Omega{u}_{m}^i=\frac{\varepsilon nc_0}{\varepsilon^\varrho}\int_\Omega\partial_r(\nu\cdot{w}_{m}^i)=-\varepsilon nc_0\int_{\partial\Omega}\nu\cdot{w}_{m}^i=\varepsilon nc_0\int_{\partial\Omega}\nu\cdot{w}_{m}^b\ .\label{eq: mean direction half space}
\end{equation}
Finally we obtain from (\ref{eq: mass conservation boundary}), 
\begin{align}
   \partial_t\left(\varepsilon^{2}\int_\Omega{u}_{m}^i-\varepsilon^{2\mu-\varrho+3}\int_{\Omega_{2\delta}}{u}^b_{m} \right)=0\ .\label{eq: only u}
\end{align}
Choosing $\varrho-2\mu=1$ and for $m=0$ we match the interior and boundary solutions
\[
\partial_t\left(\int_\Omega{u}_{0}^i-\int_{\Omega_{2\delta}}{u}^b_{0} \right)=0\ .
\]
Let us consider the mass conservation equation at the boundary, given by equation \eqref{eq: mass conservation boundary}, which we can rewrite as
\begin{equation}
\varepsilon^{2\mu+2}\int_{\Omega_{2\delta}}\partial_tu_{m}^b+c_0n\int_{\Omega_{2\delta}}\partial_r(\nu\cdot w_m^b)=0\ .\label{eq: conservation at the boundary 2}
\end{equation}
Now assume the expansions  
\begin{equation}\label{eq: expansions u and w}
w_m^b= w_0^b+\varepsilon^{2\mu+2} w_1^b+\textnormal{l.o.t}\quad \textnormal{and}\quad u_m^b=u_0^b+\varepsilon u_1^b+\textnormal{l.o.t}\ .\end{equation} 
From this, it follows that, at leading order, $w_0^b=0$ and the conservation equation becomes
\[
\int_{\Omega_{2\delta}}\partial_tu_{0}^b+c_0n\int_{\Omega_{2\delta}}\partial_r(\nu\cdot w_1^b)=0\ ,
\]
which is consistent with the leading order behavior obtained in the interior solution as stated in Proposition \ref{prop: 1}.

We can now express $f_m^b$ as
\[
f_m^b=u_m^b+c_0n\mathbf{v}\cdot w_m^b=u_0^b+\varepsilon u_1^b+\varepsilon^{2\mu+2}c_0n\mathbf{v}\cdot w_1^b+\mathcal{O}(\min\{\varepsilon^2,\varepsilon^{2\mu+3}\})\ .
\]
To verify that $f_m=f^i_m+f_m^b$ satisfies the boundary conditions in \eqref{eq: original system} we substitute the expansion and obtain
\begin{equation}
    f_m^b(0,\xi,\mathbf{v},t)-f_m^b(0,\xi,R_\xi(\mathbf{v}),t)=f_m^i(\xi,R_\xi(\mathbf{v}),t)-f_m^i(\xi,\mathbf{v},t)\ .\label{eq: boundary codit expansion}
\end{equation}
We conclude that the leading-order terms in $f_m^b$ satisfy the transport equation \eqref{eq: m-expression}, the boundary condition \eqref{eq: boundary codit expansion} and the decay condition at infinity
\[
0=\lim_{r\to\infty}f_m^b(r,\xi,\mathbf{v},t)\ .
\]

\subsection{Macroscopic equation at the boundary}

We investigate the leading order equations arising from the boundary approximation \eqref{eq: m-expression}. We proceed by substituting ${f}_m^b=\frac{1}{|S|}({u}_m^b+n\mathbf{v}\cdot{w}_m^b)$ and multiplying the whole expression by $\mathbf{v}$. Integrating with respect  to $S$ we obtain
\begin{align}
    c_0\nu & \cdot\partial_r{u}_{m}^b =-\int_S\mathbf{v}(\mathds{1}-T)\Bigl[B_{\frac{2-\alpha}{\alpha-1}}+C_{\frac{2-\alpha}{\alpha-1}} \Bigr]({u}_{m}^b+n\mathbf{v}\cdot{w}_{m}^b)\diff\mathbf{v}\nonumber\\ &=-\frac{\alpha-1}{\tau_0}n{w}_{m}^b(1-\nu_1)+\frac{\pi\tau_0^2(1-\alpha)^2}{\sin(\pi\alpha)\Gamma(\alpha)}\frac{(n^2\nu_1-|S|)}{|S|}(c_0\nu\cdot\partial_r)^{\alpha-1}{u}_{m}^b+W_{m}^b\ ,\label{eq: expansion leading order}
\end{align}
where
\[
W^b_{m}=\frac{1}{|S|}\int_S(\mathds{1}-T)n\tau_0^{\alpha-2}(\alpha-1)^2\Gamma(-\alpha+1)\mathbf{v}^2(-c_0\mathbf{v}\cdot\nu\partial_r)^{\alpha-1}{w}_{m}^b\diff\mathbf{v}\ .
\]
Substituting the expressions for $u_m^b$ and $w_m^b$ in \eqref{eq: expansions u and w} into \eqref{eq: expansion leading order} we obtain
\begin{equation}
{\varepsilon^{2\mu+2}}w_1^b=\frac{-1}{C_\alpha}\nu\cdot\partial_r u_0^b+\frac{D_\alpha}{C_\alpha}(c_0\nu\cdot\partial_r)^{\alpha-1}u_0^b\ , \label{eq: mean direction boundary}
\end{equation}
where
\[
D_\alpha=\frac{\pi\tau_0^2(1-\alpha)^2}{\sin(\pi\alpha)\Gamma(\alpha)}\frac{(n^2\nu_1-|S|)}{|S|}\ ,\  C_\alpha=\frac{\alpha-1}{\tau_0}n(1-\nu_1)\ ,
\]
and $W_0^b=0$ since, as we saw before, $w_0^b=0$.
Going back to the conservation equation at the boundary given by \eqref{eq: conservation at the boundary} we can substitute \eqref{eq: mean direction boundary} to obtain
\[
\int_{\Omega_{2\delta}}\partial_t u_0^b-{\frac{1}{\varepsilon^{2\mu+2}}}\int_{\Omega_{2\delta}}\partial_r\Bigl[\nu\cdot\Bigl( \frac{1}{C_\alpha}\nu\cdot\partial_r u_0^b-\frac{D_\alpha}{C_\alpha}(c_0\nu\cdot\partial_r)^{\alpha-1}u_0^b\Bigr)\Bigr]=0\ .
\]
\begin{rem}
For $\mu=0$ this is analogous to the classical diffusion case $\varepsilon\partial_tu-\Delta u=0$ with typical boundary layer $r=x/\sqrt{\varepsilon}$, where  $\Delta u\sim\varepsilon^{-1}\partial_r^2u_0^b$,  $u(x)=u^b(r,\xi)=u^b_0+\varepsilon u^b_1+\varepsilon^2u^b_2+\mathcal{O}(\varepsilon^3)$ and we get 
$
\varepsilon\partial_tu^b_0-\varepsilon^{-1}\partial_r^2u^b_0=0\ .
$
\end{rem}

\section{Curved boundary approximation}
In this section, we will examine the construction of a boundary layer solution, denoted again as $f_m^b$ for \eqref{eq: last}, which becomes negligible beyond a few mean free paths from the physical boundary $\partial\Omega$. In the vicinity of $\partial\Omega$, we impose that $f_m^b$
exhibits slow variation along the tangential direction to $\partial\Omega$, while allowing more rapid changes in the normal direction. To express these characteristics mathematically, we introduce the curvilinear coordinate system $(\xi,d)$ near $\partial\Omega$ defined by 
\begin{equation}
    \mathbf{x}(\xi,d)=\mathbf{x}_b(\xi)+d(\mathbf{x})\nu(\xi)\ .\label{eq: curvilinear coord}
\end{equation}

\begin{defn}[Curved boundary layer]Let $\Omega_{2\delta}=\left\{d(x)<2\delta \right\}$ be a boundary strip such that  $\mathbf{x}$ is near $\partial \Omega$, $\mathbf{x}_b(\xi)\in\partial \Omega$ and $\nu(\xi)$ is the inner unit normal at $\xi$. We consider again $r=d/\varepsilon^\varrho$.
\end{defn}
Consider now the asymptotic expansion close to the boundary as 
\begin{equation}
{f}^b=\sum_{m=0}^N\varepsilon^m{f}_m^b\left(r,\xi,\mathbf{v},t \right)\ ,
\end{equation}
and the gradient is composed of normal and tangential components, i.e.,
\begin{align}
     \nabla=\frac{\nu}{\varepsilon^\varrho}\partial_r+\frac{1}{1-\kappa(\xi)d}\mathbf{t}(\xi)\, \partial_\xi\ .
\end{align}
This expression for the gradient is coming from the following considerations. We assume $\mathbf{t}(\xi)=\partial\mathbf{x}/\partial\xi$ is the unit tangent vector and $\nu(\xi)=\partial\mathbf{x}/\partial d$ is the unit normal vector. Differentiating \eqref{eq: curvilinear coord} we have
\[
\frac{\partial\mathbf{x}}{\partial\xi}=\mathbf{t}+d\frac{\partial\nu}{\partial\xi}=(1-\kappa d)\mathbf{t}\ ,
\]
where $\partial\nu/\partial\xi=-\kappa(\xi)\mathbf{t}$ and $\kappa(\xi)$ is the curvature.

\begin{pro}\label{pro: 3.5}
There exist smooth functions ${f}_m^b:\ \mathds{R}^2_+\times\partial\Omega\times S\rightarrow\mathds{R}$ which are uniquely determined by the following system of boundary value problems:
\begin{align}
c_0\mathbf{v}\cdot\nu\partial_r{f}_{m}^b+(\mathds{1}-T)\Bigl({B}_{\frac{2-\alpha}{\alpha-1}}+{C}_{\frac{2-\alpha}{\alpha-1}}\Bigr){f}_{m}^b&= g_{m}^b\thinspace, \ \mathrm{on}\ \mathds{R}_+\ ,\label{eq: boundary value problem}\\
{f}_m^b(0,\xi,\mathbf{v},t) &={f}_m^b(0,\xi,\textnormal{R}_\xi(\mathbf{v}),t)\ .
\end{align}
Here, 
\begin{align}
    g_{m}^b=-\varepsilon^{2\mu+2}\partial_t{f}_{m}^b-\varepsilon^{2\mu+1}c_0\frac{\mathbf{v}\cdot\mathbf{t}}{1-\kappa d}\partial_\xi f_{m}^b-\varepsilon(\mathds{1}-T) {A}_{2\mu-1}{f}_{m}^b\ . \label{eq: curved boundary equation}
\end{align}

\end{pro}
Following the same procedure as before, we substitute ${f}_m^b=\frac{1}{|S|}({u}_m^b+n\mathbf{v}\cdot{w}_m^b)$ into (\ref{eq: boundary value problem}) and integrate in $S$ to obtain
\begin{equation}
\partial_r(\nu\cdot{w}_{m}^b)=\frac{1}{nc_0}\int_Sg_{m}^bd\mathbf{v}\ . \label{eq: E0 component}
\end{equation}

Knowing that
\[
\int_{\partial\Omega}c_0\nu\cdot{w}_{m}^b=\int_{\Omega_{2\delta}}c_0\nabla\cdot{w}^b_{m}\ ,
\]
from the Stokes' theorem, and considering only the normal component of the flux at the boundary, we get 
\begin{align} \int_{\partial_\Omega}\nu\cdot{w}_{m}^b & = \int_{\Omega_{2\delta}}\nabla\cdot(\nu(\nu\cdot{w}_{m}^b))\nonumber\\ &= \int_{\Omega_{2\delta}}[(\nabla\cdot\nu)(\nu\cdot{w}_{m}^b)+\nu\cdot\nabla(\nu\cdot{w}_{m}^b)]\nonumber\\ & =\int_{\Omega_{2\delta}}\left[\Delta d\, (\nu\cdot{w}_{m}^b)+\varepsilon^{-\varrho}\partial_r(\nu\cdot{w}_{m}^b)\right]\ .\label{eq: conservation curved}
\end{align}
Substituting (\ref{eq: E0 component}) into the right hand side of (\ref{eq: conservation curved}) we finally obtain
\begin{align}
\int_{\partial\Omega}\nu\cdot{w}_{m}^b&=\int_{\Omega_{2\delta}}\Bigl[\Delta d\, (\nu\cdot w_{m}^b)+\frac{\varepsilon^{-\varrho}}{nc_0}\int_Sg_{m}^bd\mathbf{v}\Bigr]\nonumber\\ &=\int_{\Omega_{2\delta}}\Delta d\, (\nu\cdot w_{m}^b)-\frac{\varepsilon^{-\varrho+2\mu+2}}{nc_0}\, \int_{\Omega_{2\delta}}\partial_t{u}_{m}^b\ ,\label{eq: mean direction curved}
\end{align}
where the last equality comes from \eqref{eq: curved boundary equation} where we have
{\begin{align}
    \int_Sg_{m}^bd\mathbf{v}=-\varepsilon^{2\mu+2}\partial_t\left({u}_{m}^b+\textnormal{l.o.t.} \right)\ .
\end{align}
Finally, from (\ref{eq: mean direction curved}) and (\ref{eq: mean direction half space}) we obtain an equation satisfied by the interior and the boundary solutions, i.e.,
\begin{equation}
 \partial_t\left(\varepsilon\int_\Omega{u}_{m}^i+{\varepsilon^{-\varrho+2+2\mu}}\int_{\Omega_{2\delta}}{u}^b_{m} \right)=\int_{\Omega_{2\delta}}\Delta d\, (\nu\cdot w_{m}^b)\ .\label{eq: conservation of partial}
\end{equation}
\iffalse
\begin{pro}
Given that the solution of the system (\ref{eq: original system}) in the interior of the domain, ${u}_m^i$, and the solution near the boundary ${u}_m^b$ satisfy the relation (\ref{eq: conservation of partial}), then it follows the conservation of mass property:
\[
\int_\Omega\int_S({f}_{(N)}^i+{f}^b_{(N)})(\mathbf{x},t,\mathbf{v})d\mathbf{v} dx=\int_\Omega\int_S\delta_0(\mathbf{x})d\mathbf{v} dx.
\]
\end{pro}
\fi

As in Section \ref{sec: boundary}, we choose $\varrho-2\mu=1$.

Considering again \eqref{eq: conservation curved} and \eqref{eq: mean direction curved}, we can write a conservation equation at the boundary as
\begin{equation}\label{eq: conservation curved boundary}
\varepsilon^{2\mu+2}\int_{\Omega_{2\delta}}\partial_t u_m^b+\int_{\Omega_{2\delta}}\partial_r\nu\cdot w_m^b=0\ ,
\end{equation}
and this leads to $u_m^b=u_0^b+\varepsilon u_1^b+l.o.t$ and $w_m^b=w_0^b+\varepsilon^{2\mu+2}w_1^b+l.o.t$ and therefore we write
\[
\int_{\Omega_{2\delta}}\partial_t u_0^b+\int_{\Omega_{2\delta}}\partial_r\nu\cdot w_1^b=0\ .
\]

Finally, we conclude that $f_m^b$ is given by \eqref{eq: boundary codit expansion} and it satisfies the boundary condition in Proposition \ref{pro: 3.5}. Then we conclude that the leading order terms in $f_m^b$ satisfy the transport equation \eqref{eq: boundary value problem}, the boundary conditions and the condition
\[
0=\lim_{r\to\infty}f^b_m(r,\xi,\mathbf{v},t)\ .
\]

\section{Reflective boundary: A case of study}

\subsection{Boundary solution}
The boundary correction \(f^b\) satisfies \eqref{eq: boundary value problem} irrespective of the boundary condition that has been imposed. In the following we are going to discuss a notion of solution for this equation following \cite{larsen}.

Consider the auxiliary problem in the half-space \(r \in (0,\infty)\),
\begin{equation}
\label{eq:x}
c_0 \,\bv\cdot \nu\partial_r\,f^b_{m}
+(\mathds{1} - T)\Bigl({B}_{\frac{2-\alpha}{\alpha-1}}+{C}_{\frac{2-\alpha}{\alpha-1}}\Bigr)f^b_{m}
=g^b_{m}\ ,
\end{equation}
subject to
\[
\lim_{r\to\infty}\bigl|f^b_m\bigr|<\infty,
\quad
f^b_m\bigl(0,{\xi},\bv,t\bigr)=\ell\bigl(\xi,\bv,t\bigr).
\]
Note that \(f^b_{m} = 1\) is an exact solution to the first equation.
We shall assume that the solution to \eqref{eq:x} takes the form
\begin{equation}
\label{eq:starstar}
f^b\bigl(r,\xi,\bv,t\bigr)
=
\int_{\bv'\cdot \nu>0}
W\bigl(\xi,\bv',t\bigr)\,
\ell\bigl(\xi,\bv',t\bigr)\,\diff\bv'
+
\int_{\bv'\cdot \nu>0}
G\bigl(r,\xi,\bv',\bv,t\bigr)\,
\ell\bigl(\xi,\bv',t\bigr)\,\diff\bv'\,,
\end{equation}
where 
\[
G\bigl(r,\xi,\bv',\bv,t\bigr)
\longrightarrow 0
\quad\text{as}\quad r\to\infty.
\]
This means that the solution operator consists of the leading behavior for \(r\to\infty\) and a decaying remainder. If \(W\) and \(G\) exist, they must satisfy
\begin{equation}
1 = \int_{\bv'\cdot \nu>0}
W\bigl(\xi,\bv',t\bigr)\,\diff\bv',\label{eq: first condition}
\end{equation}
\begin{equation}
0 = \int_{\bv'\cdot \nu>0}
G\bigl(r,\xi,\bv',\bv,t\bigr)\,\diff\bv'.
\label{eq: second condition}
\end{equation}
Such \(W,\,G\) are explicitly known for some problems \cite{bell1970nuclear,davison1957neutron}. 
A solution \(f^b\) of the form \eqref{eq:starstar} tends to \(0\) as \(r\to\infty\), provided that
\begin{equation}
\int_{\bv'\cdot\nu>0}
W\bigl(\xi,\bv',t\bigr)\,
\ell\bigl(\xi,\bv',t\bigr)\,\diff\bv'
=0\ . \label{eq: necessary condition W}
\end{equation}
We now derive some properties of the reflection operator $\textnormal{R}: X^+ \longrightarrow X^-$,
where \(X^\pm\) consists of functions \(f(\xi,\bv,t)\) defined for \(\xi\) on the boundary, \(t>0\), and \(\bv\cdot \nu\ge 0\).  We define
\begin{equation}
(\textnormal{R}\,\ell)\bigl(\xi,\bv,t\bigr)
=
\int_{\bv'\cdot \nu>0}
W\bigl(\xi,\bv',t\bigr)\,\ell\bigl(\xi,\bv',t\bigr)\,\diff\bv'
+
\int_{\bv'\cdot \nu>0}
G\bigl(0,\xi,\bv',\bv,t\bigr)\,\ell\bigl(\xi,\bv',t\bigr)\,\diff\bv',\label{eq: reflection}
\end{equation}
valid for \(\bv\cdot\nu<0\).

\medskip

By integrating 
\eqref{eq:x} over \(\bv\), we obtain an important property of \(\textnormal{R}\).  Namely,
\[
0
=
\partial_r\,\int \bigl(\bv\cdot \nu\bigr)\,
f_m^b\bigl(r,\xi,\bv,t\bigr)\,\diff\bv,
\]
so the integral on the right hand side is independent of $r$,
\begin{equation}\label{eq: rand expression}
\int \bigl(\bv\cdot \nu\bigr)\,
f_m^b\bigl(0,\xi,\bv,t\bigr)\,d\bv=\int(\bv\cdot\nu)f_m^b(r,\xi,\bv,t)d\bv\ .
\end{equation}
As \(r\to \infty\), \(f_m^b\) tends to a constant, so the integral in \eqref{eq: rand expression} tends to \(0\). Thus
\[
\int \bigl(\bv\cdot \nu\bigr)\,
f_m^b\bigl(0,\xi,\bv,t\bigr)\,d\bv
=0.
\]
Breaking the integral into  \(\bv\cdot\nu>0\) and \(\bv\cdot\nu<0\), we obtain
\[
0
=
\int_{\bv\cdot\nu>0}
\bigl(\bv\cdot \nu\bigr)\,
f_m^b\bigl(0\xi,\bv,t\bigr)\,\diff\bv
+
\int_{\bv\cdot\nu<0}
\bigl(\bv\cdot \nu\bigr)\,
f_m^b\bigl(0,\xi,\bv,t\bigr)\,\diff\bv\ ,
\]
and identifying
\(\ell(\xi,\bv,t)\) for \(\bv\cdot\nu>0\)
and
\((\textnormal{R}\ell)(\xi,\bv,t)\) for \(\bv\cdot\nu<0\), 
this yields, 
\begin{equation}
\int_{\bv\cdot\nu>0}(\bv\cdot\nu)\ell(\xi,\bv,t)\diff\bv=\int_{\bv\cdot\nu<0}|\bv\cdot\nu|(\textnormal{R}\ell)(\xi,\bv,t)\diff\bv\ .\label{eq: 4.25}
\end{equation}

The adjoint operator \(\widetilde{\textnormal{R}}: (X^-)^*\to (X^+)^*\) thus satisfies
\[
\widetilde{\textnormal{R}}\bigl(|\bv\cdot \nu|\bigr)
=
\bv\cdot \nu,
\quad
\bv\cdot \nu>0.
\]
Similarly, from  conditions \eqref{eq: first condition}, \eqref{eq: second condition} and \eqref{eq: reflection} we have
\[
\textnormal{R}(1) = 1
\quad
\text{for} \, \bv\cdot\nu<0.
\]

In the next step we solve a half-space problems for \(f_m^b\) in terms of the incoming density
\[
\ell(\xi,\bv,t) = f_m^b\bigl(0,\xi,\bv,t\bigr),
\quad
\bv\cdot\nu>0.
\]
The above formalism leads to necessary conditions for the density \(\ell\), in particular that it satisfies the condition \eqref{eq: necessary condition W} above.  These conditions determine \(\ell\) uniquely and lead to a boundary condition for \(u_0\), the solution of the macroscopic diffusion equation. 
%\textcolor{red}{(We should think about exterior conditions, too!)}

\subsection{Reflective boundary condition}
We consider the \emph{prompt reflection operator} $\textnormal{P}: X^- \longrightarrow X^+$,
written as 
\begin{equation}
\label{eq:P-def}
\bigl(\bv\cdot \nu\bigr)\,
\bigl(\textnormal{P}f_m^b\bigr)\bigl(\xi,\bv,t\bigr)
=
\int_{\bv'\cdot \nu<0}
\bigl|\bv'\cdot \nu\bigr|\,
p\bigl(\xi,\bv',\bv,t\bigr)
f_m^b\bigl(\xi,\bv',t\bigr)\,\diff\bv',
\end{equation}
for \(\xi\) on the boundary and \(\bv\cdot \nu>0\). The kernel \(p\) is nonnegative and can be for instance a Dirac-\(\delta\) in the case of specular reflection.

Conservation of particles under prompt reflection is assured if
\[
1
=
\int_{\bv'\cdot\nu>0}
p\bigl(\xi,\bv',\bv,t\bigr)\,
\diff\bv',
\quad
\bv'\cdot \nu>0.
\]
Using this property, \(\textnormal{P}\) satisfies
\begin{equation}
\int_{\bv\cdot \nu>0}
\bigl(\bv\cdot \nu\bigr)\,
\bigl(\textnormal{P}f_m^b\bigr)(\xi,\bv,t)\,\diff\bv
=
\int_{\bv'\cdot \nu<0}
\bigl|\bv'\cdot \nu\bigr|\,
f_m^b\bigl(\xi,\bv',t\bigr)\,\diff\bv',\label{eq: condition P}
\end{equation}
and its adjoint \(\widetilde{\textnormal{P}}\) acts by
$\widetilde{\textnormal{P}}\bigl(\bv\cdot\nu\bigr)
=
\bigl|\bv\cdot \nu\bigr|,$ for
$\bv\cdot\nu<0$.

%\textcolor{red}{(We might later include small \emph{delayed} reflections, corresponding to \(P_1\) and \(Q\) in Larsen.)}

To determine the boundary conditions for \(u_0\), we first consider the problem \eqref{eq: boundary value problem} for \(f^b_m\): 
The boundary condition is written as
\[
f^b_m\bigl(0,\xi,\bv,t\bigr)
+
f^i_m\bigl(\xi,t,\bv\bigr)
=
f^b_m\bigl(0,\xi,\textnormal{R}_\xi(\bv),t\bigr)
+
f^i_m\bigl(\xi,t,\textnormal{R}_\xi(\bv)\bigr).
\]
With $f^i_m\bigl(\xi,t,\bv\bigr) = u_0\bigl(\xi,t,\bv\bigr)$, $f^i_m\bigl(\xi,t,\textnormal{R}_\xi(\bv)\bigr) = \textnormal{PR}\,u_0\bigl(\xi,t,\bv\bigr)$ and $f^b_m\bigl(0,\xi,\textnormal{R}_\xi(\bv),t\bigr) = \textnormal{PR} f^b_m\bigl(0,\xi,\bv,t\bigr)$, for the specular reflection we find $f^b_m - \textnormal{PR}\,f^b_m=\textnormal{PR}\,u_0-u_0$.
We write this as an integral equation for the unknown  
\[
\ell_0\bigl(\xi,\bv,t\bigr)
=
f^b_m\bigl(0,\xi,\bv,t\bigr),
\]
which involves the operator $\mathds{1} - \textnormal{PR}$:
\begin{equation} \label{PRequation}
\bigl(\mathds{1} - \textnormal{PR}\bigr)\bigl(\ell_0 + u_0\bigr)\bigl(\xi,\bv,t\bigr)
=
0.
\end{equation}\\
Recall from our earlier discussion that therefore \eqref{eq: boundary value problem} admits a unique solution in the half-space if and only if the condition \eqref{eq: necessary condition W} holds for \(\ell = \ell_0\). 

We further note that equations \eqref{eq: 4.25}  and \eqref{eq: condition P} imply for any $\ell_0$ that
\begin{equation}\label{solvecond}
0
=
\int_{\bv\cdot\nu>0}
\bigl(\bv\cdot \nu\bigr)\,\bigl(\mathds{1} - \textnormal{PR}\bigr)\bigl(\ell_0 + u_0\bigr)\,\diff\bv.
\end{equation}
\iffalse
\textcolor{red}{\[
\int(\bv\cdot\nu)(Pf_m^b)d\bv-\int(\bv\cdot\nu)\ell d\bv=0
\]}
\fi
This shows that the function $\Theta^*=\bv\cdot\nu$ is in the null space of the adjoint operator $\widetilde{\mathds{1}-\textnormal{PR}}$ of $\mathds{1}-\textnormal{PR}$. Therefore the operator $\mathds{1}-\textnormal{PR}$ also has a null space, and we denote by $\Theta=\Theta(\xi,\bv,t)\in X^+$ a nonzero function in it. In the following we assume that, in fact, $\Theta>0$ is normalised as
\begin{equation}\label{Wnormalize}
1 
=
\int_{\bv\cdot\nu>0}
W\bigl(\xi,\bv,t\bigr)\,\Theta\bigl(\xi,\bv,t\bigr)\,\diff\bv\ ,
\end{equation} and that it spans the null space of $\mathds{1}-\textnormal{PR}$. We further assume that $\lambda=0$ is an isolated eigenvalue.

Because the null space is spanned by $\Theta$, the general solution of \eqref{PRequation} is of the form
\begin{equation}\label{l0form}
\ell_0
=
a_0\,\Theta 
-
u_0, \quad \textnormal{with}\quad a_0=a_0(\xi,t),
\quad
\bv\cdot \nu(\xi)>0,
\end{equation}
for some $a_0$. We determine $a_0$ with the help of \eqref{eq: necessary condition W} and \eqref{Wnormalize}:
\[
0
=
\int W\,\ell_0 \diff\bv'
=
a_0 \int W\,\Theta \diff\bv'
-
u_0 \int W \diff\bv'
=
a_0
-
u_0.
\]
This shows $a_0=u_0$, and therefore \eqref{l0form} becomes
\[
\ell_0\bigl(\xi,\bv,t\bigr)
=
u_0\bigl(\xi,t\bigr)\,\Bigl[\Theta\bigl(\xi,\bv,t\bigr) - 1\Bigr],
\quad
\bv\cdot \nu(\xi)>0.
\]
This determines \(\ell_0\) up to the multiplicative factor \(u_0\).

\smallskip

To determine \(u_0\), we consider the problem for \(f_{m+1}\).  The boundary condition satisfied by \(\ell_1 = f^b_{m+1}(\xi,0,\bv,t)\) is
\begin{equation}\label{eq: expresion ell 1}
\bigl(\mathds{1} - \textnormal{PR}\bigr)\bigl(\ell_1 + u_1\bigr)
=h_1+{(\textnormal{P}-\mathds{1})\bv\cdot w_1^b},
\quad
\bv\cdot\nu>0,
\end{equation}
for a given  $h_1$. We multiply it by $\bv\cdot\nu$ and integrate over the half-space $\bv\cdot\nu >0$. As in equation \eqref{solvecond}, we find that the left hand side of \eqref{eq: expresion ell 1} vanishes and we have the following solvability condition for the data $h_1$,
\[
\int_{\bv\cdot\nu>0}(\bv\cdot\nu)(\mathds{1}-\textnormal{P})(\bv\cdot w_1^i)\,\diff \bv=\int_{\bv\cdot\nu>0}
\bigl(\bv\cdot \nu\bigr)\,h_1\,\diff\bv
\]
then, from \eqref{eq: mean direction} we write
\begin{align}
\int_{\bv\cdot\nu>0}(\bv\cdot\nu)(\mathds{1}-\textnormal{P})(\bv\cdot w_1^i)\,\diff \bv &= \mathcal{H}\, u_0^i\ ,
\end{align}
where
\[
\mathcal{H}\, u_0^i\coloneqq\Bigl(\int_{\mathbf{v}\cdot\nu>0}(\bv\cdot\nu)(\mathds{1}-P)[\chi\, \bv\cdot\nabla\rho-C_\alpha\bv\cdot\nabla^{\alpha-1}]\diff \bv \,\Bigr)u_0^i\ ,
\]
$C_\alpha$ and $\chi$ are defined in \eqref{eq: constants}.
Therefore
\[
\mathcal{H}\, u_0^i=\nu\cdot\int_{\bv\cdot\nu>0}\bv h_1\diff \bv \ .
\]
This provides the relevant boundary condition for the interior problem \eqref{eq: main at intro}, assuring conservation of particles at the boundary.

\section{Conclusions}
%\begin{itemize}
%    \item Think about exterior boundary conditions. Modify (4.17) in Larsen's.
%    \item Comment on the choice of fractional radial derivative. Whether we allow particles to leave the domain and come back, or if we kill particles outside the domain etc.
%\end{itemize}
Employing boundary layer techniques as in \cite{alt1980singular, larsen}, we systematically constructed corrections near the domain boundary, accounting for the singular behavior introduced by the nonlocal nature of fractional operators. This analysis revealed the persistence of boundary influence even far from the boundary—a phenomenon unique to nonlocal models. Through the matching procedure, we demonstrated that the interior and boundary solutions are consistent at leading order. This matching ensures global mass conservation and validates the asymptotic expansions within both the bulk and the boundary layers.

We explored the reflective boundary condition in detail, showing how classical notions (such as specular reflection) must be reformulated when dealing with fractional derivatives. %The use of operators like $\textnormal{R}$ and $\textnormal{P}$, and the analysis of their adjoints. 
Our analysis consistently recovers the classical Patlak–Keller–Segel equation and Neumann boundary conditions in the limit as the fractional exponent $\alpha\to 2$, offering a smooth transition between the classical and fractional models.

\bibliographystyle{plain}
\bibliography{boundary_conditions}

\begin{thebibliography}{10}

\bibitem{alt1980singular}
Wolfgang Alt.
\newblock Singular perturbation of differential integral equations describing
  biased random walks.
\newblock {\em Journal f\"{u}r die reine und angewandte {M}athematik},
  322:15--41, 1981.

\bibitem{baeumer2018reprint}
Boris Baeumer, Mih{\'a}ly Kov{\'a}cs, Mark~M Meerschaert, and Harish
  Sankaranarayanan.
\newblock Boundary conditions for fractional diffusion.
\newblock {\em Journal of Computational and Applied Mathematics}, 339:414--430,
  2018.

\bibitem{baeumer2018fractional}
Boris Baeumer, Mih{\'a}ly Kov{\'a}cs, and Harish Sankaranarayanan.
\newblock Fractional partial differential equations with boundary conditions.
\newblock {\em Journal of Differential Equations}, 264(2):1377--1410, 2018.

\bibitem{bell1970nuclear}
George~I Bell and Samuel Glasstone.
\newblock Nuclear reactor theory.
\newblock Technical report, US Atomic Energy Commission, Washington, DC (United
  States), 1970.

\bibitem{davison1957neutron}
Boris Davison.
\newblock {\em Neutron transport theory}.
\newblock Oxford: Clarendon Press, 1957.

\bibitem{robots}
Gissell Estrada-Rodriguez and Heiko Gimperlein.
\newblock Interacting particles with l{\'e}vy strategies: limits of transport
  equations for swarm robotic systems.
\newblock {\em SIAM Journal on Applied Mathematics}, 80(1):476--498, 2020.

\bibitem{pks}
Gissell Estrada-Rodriguez, Heiko Gimperlein, and Kevin~J Painter.
\newblock Fractional {P}atlak--{K}eller--{S}egel equations for chemotactic
  superdiffusion.
\newblock {\em SIAM Journal on Applied Mathematics}, 78(2):1155--1173, 2018.

\bibitem{keller1970initiation}
Evelyn~F Keller and Lee~A Segel.
\newblock Initiation of slime mold aggregation viewed as an instability.
\newblock {\em Journal of theoretical biology}, 26(3):399--415, 1970.

\bibitem{korobkova2004molecular}
Ekaterina Korobkova, Thierry Emonet, Jose~MG Vilar, Thomas~S Shimizu, and
  Philippe Cluzel.
\newblock From molecular noise to behavioural variability in a single
  bacterium.
\newblock {\em Nature}, 428(6982):574--578, 2004.

\bibitem{larsen}
Edward~W Larsen.
\newblock Asymptotic theory of the linear transport equation for small mean
  free paths. ii.
\newblock {\em SIAM Journal on Applied Mathematics}, 33(3):427--445, 1977.

\bibitem{li2008persistent}
Liang Li, Simon~F N{\o}rrelykke, and Edward~C Cox.
\newblock Persistent cell motion in the absence of external signals: a search
  strategy for eukaryotic cells.
\newblock {\em PLoS {O}ne}, 3(5):e2093, 2008.

\bibitem{you2020asymptotically}
Huaiqian You, XinYang Lu, Nathaniel Task, and Yue Yu.
\newblock An asymptotically compatible approach for neumann-type boundary
  condition on nonlocal problems.
\newblock {\em ESAIM: Mathematical Modelling and Numerical Analysis},
  54(4):1373--1413, 2020.

\end{thebibliography}

\end{document}